\documentclass{tran-l}

\numberwithin{equation}{section}
\newtheorem{thm}{Theorem}[section]

\newtheorem{prop}[thm]{Proposition}

\theoremstyle{definition}
\newtheorem{defn}{Definition}[section]
\newtheorem{rem}{Remark}[section]

\theoremstyle{remark}

\newcommand{\X}{\mathfrak{X}}
\newcommand{\s}{\mathfrak{S}}
\newcommand{\g}{\mathfrak{g}}
\newcommand{\W}{\mathcal{W}}

\newfont{\w}{msbm9 scaled\magstep1}
\def\R{\mbox{\w R}}
\newcommand{\norm}[1]{\left\Vert#1\right\Vert ^2}
\newcommand{\nJ}{\norm{\nabla J}}
\newcommand{\nN}{\norm{N}}
\newcommand{\ad}{{\rm ad}}

\newcommand{\thmref}[1]{Theorem~\ref{#1}}
\newcommand{\propref}[1]{Proposition~\ref{#1}}

\newcommand{\dfnref}[1]{Definition~\ref{#1}}

\begin{document}

\title[A Lie group as a quasi-K\"ahler manifold with Norden metric]
{A Lie group as a 4-dimensional quasi-K\"ahler manifold with
Norden metric}

\author{Kostadin Gribachev, Mancho Manev, Dimitar Mekerov}

\address{University of Plovdiv, Faculty of Mathematics and
Informatics, Department of Geometry, 236 Bulgaria blvd., Plovdiv
4003, BULGARIA}

\email{costas@pu.acad.bg, mmanev@yahoo.com, mircho@pu.acad.bg}



\begin{abstract}
A 4-parametric family of 4-dimensional quasi-K\"ahler manifolds
with Norden metric is constructed on a Lie group. This family is
characterized geometrically. The condition for such a 4-manifold
to be isotropic K\"ahler is given.
\end{abstract}

\maketitle

\section*{Introduction}

It is a fundamental fact that on an almost complex manifold with
Hermitian metric (almost Hermitian manifold), the action of the
almost complex structure on the tangent space at each point of the
manifold is isometry. There is another kind of metric, called a
Norden metric or a $B$-metric on an almost complex manifold, such
that action of the almost complex structure is anti-isometry with
respect to the metric. Such a manifold is called an almost complex
manifold with Norden metric \cite{GaBo} or with $B$-metric
\cite{GaGrMi}. See also \cite{GrMeDj} for generalized
$B$-manifolds. It is known \cite{GaBo} that these manifolds are
classified into eight classes.

The purpose of the present paper is to exhibit, by construction,
almost complex structures with Norden metric on Lie groups as
4-manifolds, which are of certain classes, called quasi-K\"ahler
manifold with Norden metrics. It is proved that the constructed
4-manifold is isotropic K\"ahler (a notion introduced by
Garc\'ia-R\'io and Matsushita \cite{GRMa}) if and only if it is
scalar flat.


\section{Almost Complex Manifolds with Norden Metric}\label{sec_1}

\subsection{Preliminaries}\label{sec-prelim}

Let $(M,J,g)$ be a $2n$-dimensional almost complex manifold with
Norden metric, i.e. $J$ is an almost complex structure and $g$ is
a metric on $M$ such that
\begin{equation}\label{Jg}
J^2X=-X, \qquad g(JX,JY)=-g(X,Y)
\end{equation}
for all differentiable vector fields $X$, $Y$ on $M$, i.e. $X, Y
\in \X(M)$.

The associated metric $\tilde{g}$ of $g$ on $M$ given by
$\tilde{g}(X,Y)=g(X,JY)$ for all $X, Y \in \X(M)$ is a Norden
metric, too. Both metrics are necessarily of signature $(n,n)$.
The manifold $(M,J,\tilde{g})$ is an almost complex manifold with
Norden metric, too.

Further, $X$, $Y$, $Z$, $U$ ($x$, $y$, $z$, $u$, respectively)
will stand for arbitrary differentiable vector fields on $M$
(vectors in $T_pM$, $p\in M$, respectively).

The Levi-Civita connection of $g$ is denoted by $\nabla$. The
tensor filed $F$ of type $(0,3)$ on $M$ is defined by
\begin{equation}\label{F}
F(X,Y,Z)=g\bigl( \left( \nabla_X J \right)Y,Z\bigr).
\end{equation}
It has the following symmetries
\begin{equation}\label{F-prop}
F(X,Y,Z)=F(X,Z,Y)=F(X,JY,JZ).
\end{equation}

Further, let $\{e_i\}$ ($i=1,2,\dots,2n$) be an arbitrary basis of
$T_pM$ at a point $p$ of $M$. The components of the inverse matrix
of $g$ are denoted by $g^{ij}$ with respect to the basis
$\{e_i\}$.

The Lie form $\theta$ associated with $F$ is defined by
\begin{equation}\label{theta}
\theta(z)=g^{ij}F(e_i,e_j,z).
\end{equation}

A classification of the considered manifolds with respect to $F$
is given in \cite{GaBo}. Eight classes of almost complex manifolds
with Norden metric are characterized there according to the
properties of $F$. The three basic classes are given as follows
\begin{equation}\label{class}
\begin{array}{l}
\W_1: F(x,y,z)=\frac{1}{4n} \left\{
g(x,y)\theta(z)+g(x,z)\theta(y)\right. \\[4pt]
\phantom{\mathcal{W}_1: F(x,y,z)=\frac{1}{4n} }\left.
    +g(x,J y)\theta(J z)
    +g(x,J z)\theta(J y)\right\};\\[4pt]
\W_2: \mathop{\s} \limits_{x,y,z}
F(x,y,J z)=0,\quad \theta=0;\\[8pt]
\W_3: \mathop{\s} \limits_{x,y,z} F(x,y,z)=0,
\end{array}
\end{equation}
where $\s $ is the cyclic sum by three arguments.

The special class $\W_0$ of the K\"ahler manifolds with Norden
metric belonging to any other class is determined by the condition
$F=0$.

\subsection{Curvature properties}\label{sec-curv}

Let $R$ be the curvature tensor field of $\nabla$ defined by
\begin{equation}\label{R}
    R(X,Y)Z=\nabla_X \nabla_Y Z - \nabla_Y \nabla_X Z -
    \nabla_{[X,Y]}Z.
\end{equation}
The corresponding tensor field of type $(0,4)$ is
determined as follows
\begin{equation}\label{R04}
    R(X,Y,Z,U)=g(R(X,Y)Z,U).
\end{equation}
The Ricci tensor $\rho$ and the scalar curvature $\tau$ are
defined as usually by
\begin{equation}\label{rho-tau}
    \rho(y,z)=g^{ij}R(e_i,y,z,e_j),\qquad \tau=g^{ij}\rho(e_i,e_j).
\end{equation}

It is well-known that the Weyl tensor $W$ on a $2n$-dimensional
pseudo-Riemannian manifold ($n\geq 2$) is given by
\begin{equation}\label{W}
    W=R-\frac{1}{2n-2}\psi_1(\rho)-\frac{\tau}{2n-1}\pi_1,
\end{equation}
where
\begin{equation}\label{psi-pi}
    \begin{array}{l}
      \psi_1(\rho)(x,y,z,u)=g(y,z)\rho(x,u)-g(x,z)\rho(y,u) \\
      \phantom{\psi_1(\rho)(x,y,z,u)}+\rho(y,z)g(x,u)-\rho(x,z)g(y,u); \\
      \pi_1=\frac{1}{2}\psi_1(g)=g(y,z)g(x,u)-g(x,z)g(y,u). \\
    \end{array}
\end{equation}

Moreover, the Weyl tensor $W$ is zero if and only if the manifold
is conformally flat.

Let $\alpha=\{x,y\}$ be a non-degenerate 2-plane (i.e.
$\pi_1(x,y,y,x) \neq 0$) spanned by vectors $x, y \in T_pM, p\in
M$. Then, it is known, the sectional curvature of $\alpha$ is
defined by the following equation
\begin{equation}\label{k}
    k(\alpha)=k(x,y)=\frac{R(x,y,y,x)}{\pi_1(x,y,y,x)}.
\end{equation}

The basic sectional curvatures in $T_pM$ with an almost complex
structure and a Norden metric $g$ are
\begin{itemize}
    \item \emph{holomorphic sectional curvatures} if $J\alpha=\alpha$;
    \item \emph{totally real sectional curvatures} if
    $J\alpha\perp\alpha$ with respect to $g$.
\end{itemize}

In \cite{GrDjMe}, a \emph{holomorphic bisectional curvature}
$h(x,y)$ for a pair of holomorphic 2-planes $\alpha_1=\{x,Jx\}$
and $\alpha_2=\{y,Jy\}$ is defined by
\begin{equation}\label{h}
    h(x,y)=-\frac{R(x,Jx,y,Jy)}
    {\sqrt{\pi_1(x,Jx,x,Jx)\pi_1(y,Jy,y,Jy)}},
\end{equation}
where $x$, $y$ do not lie along the totally isotropic directions,
i.e. the both of the couples $\bigl(g(x,x), g(x,Jx)\bigr)$ and
$\bigl(g(y,y), g(y,Jy)\bigr)$ are different from the couple
$\left(0,0\right)$. The holomorphic bisectional curvature is
invariant with respect to the basis of the 2-planes $\alpha_1$ and
$\alpha_2$. In particular, if $\alpha_1=\alpha_2$, then the
holomorphic bisectional curvature coincides with the holomorphic
sectional curvature of the 2-plane $\alpha_1=\alpha_2$.

\subsection{Isotropic K\"ahler manifolds}\label{sec-iK}
The square norm $\nJ$ of $\nabla J$ is defined in \cite{GRMa} by
\begin{equation}\label{snorm}
    \nJ=g^{ij}g^{kl}
    g\bigl(\left(\nabla_{e_i} J\right)e_k,\left(\nabla_{e_j}
    J\right)e_l\bigr).
\end{equation}

Having in mind the definition \eqref{F} of the tensor $F$ and the
properties \eqref{F-prop}, we obtain the following equation for
the square norm of $\nabla J$
\begin{equation}\label{snormF}
    \nJ=g^{ij}g^{kl}g^{pq}F_{ikp}F_{jlq},
\end{equation}
where $F_{ikp}=F(e_i,e_k,e_p)$.

\begin{defn}[\cite{MekMan}]\label{iK}
An almost complex manifold with Norden metric satisfying the
condition $\nJ=0$ is called an \emph{isotropic K\"ahler manifold
with Norden metric}.
\end{defn}

\begin{rem}
It is clear, if a manifold belongs to the class $\W_0$, then
it is isotropic K\"ahlerian but the inverse statement is not
always true.
\end{rem}

\subsection{Quasi-K\"ahler manifolds with Norden metric}\label{sec-qK}
The only class of the three basic classes, where the almost
complex structure is not integrable, is the class $\W_3$ -- the
class of the \emph{quasi-K\"ahler manifolds with Norden metric}.

Let us remark that the definitional condition from \eqref{class}
implies the vanishing of the Lie form $\theta$ for the class
$\W_3$.

In \cite{Me} it is proved that on every $\W_3$-manifold the
curvature tensor $R$ satisfies the following identity
\begin{equation}\label{12}
\begin{array}{l}
  R(X,JZ,Y,JU)+R(X,JY,U,JZ)+R(X,JY,Z,JU)
  \\[4pt]
  +R(X,JZ,U,JY)+R(X,JU,Y,JZ)+R(X,JU,Z,JY)
  \\[4pt]
  +R(JX,Z,JY,U)+R(JX,Y,JU,Z)+R(JX,Y,JZ,U)
  \\[4pt]
  +R(JX,Z,JU,Y)+R(JX,U,JY,Z)+R(JX,U,JZ,Y)
  \\[4pt]
  =-\mathop{\s}\limits_{X,Y,Z}
    g\Bigl(\bigl(\nabla_X J\bigr)Y+\bigl(\nabla_Y J\bigr)X,
     \bigl(\nabla_Z J\bigr)U+\bigl(\nabla_U J\bigr)Z\Bigr). \\[4pt]
\end{array}
\end{equation}

As it is known, a manifold with Levi-Civita connection $\nabla$ is
\emph{locally symmetric} if and only if its curvature tensor is
parallel, i.e. $\nabla R=0$. According to \eqref{W},
\eqref{psi-pi}, $\nabla g=0$ and \eqref{12}, we establish the
truthfulness of the following
\begin{thm}\label{Prop1.1}
Every $2n$-dimensional $\W_3$-manifold with vanishing Weyl tensor
and constant scalar curvature is locally symmetric.
\end{thm}

Let us remark that a necessary condition for a $\W_3$-manifold to
be isotropic K\"ahlerian is given by the following
\begin{thm}[\cite{MekMan}]
Let $(M,J,g)$, $\dim{M}\geq 4$, be a $\W_3$-manifold and
$R(x,Jx,y,Jy)$ $=0$ for all $x,y \in T_pM$. Then $(M,J,g)$ is an
isotropic K\"ahler $\W_3$-manifold.
\end{thm}


\section{The Lie Group as a 4-Dimensional $\W_3$-Manifold}
\label{sec_2}

Let $V$ be a 4-dimensional vector space and consider the structure
of the Lie algebra defined by the brackets $
[E_i,E_j]=C_{ij}^kE_k, $ where $\{E_1,E_2,E_3,E_4\}$ is a basis of
$V$ and $C_{ij}^k\in \R$. Then the Jacobi identity for $C_{ij}^k$
\begin{equation}\label{Jac}
    C_{ij}^k C_{ks}^l+C_{js}^k C_{ki}^l+C_{si}^k C_{kj}^l=0
\end{equation}
holds.

Let $G$ be the associated connected Lie group and
$\{X_1,X_2,X_3,X_4\}$ is a global basis of left invariant vector
fields induced by the basis of $V$. Then we define an almost
complex structure by the conditions
\begin{equation}\label{J}
JX_1=X_3,\quad JX_2=X_4,\quad JX_3=-X_1,\quad JX_4=-X_2.
\end{equation}
Let us consider the left invariant metric defined by the the
following way
\begin{equation}\label{g}
\begin{array}{l}
  g(X_1,X_1)=g(X_2,X_2)=-g(X_3,X_3)=-g(X_4,X_4)=1, \\[4pt]
  g(X_i,X_j)=0\quad \text{for}\quad i\neq j. \\
\end{array}
\end{equation}
The introduced metric is a Norden metric because of \eqref{J}.

In this way, the induced 4-dimensional manifold $(G,J,g)$ is an
almost complex manifold with Norden metric. We then determine the
structural constants $C_{ij}^k$ such that the manifold $(G,J,g)$
is of class $\W_3$.

Let $\nabla$ be the Levi-Civita connection of $g$. Then the
following well-known condition is valid
\begin{equation}\label{LC}
\begin{array}{l}
    2g(\nabla_X Y,Z)=Xg(Y,Z)+Yg(X,Z)-Zg(X,Y)\\[4pt]
    \phantom{2g(\nabla_X Y,Z)=}+g([X,Y],Z)+g([Z,X],Y)+g([Z,Y],X).
\end{array}
\end{equation}

Having in mind \eqref{J}, \eqref{g} and \eqref{LC} we get
immediately the following equation for the tensor $F$ defined by
\eqref{F}
\begin{equation}\label{F_3}
\begin{array}{l}
    \mathop{\s} \limits_{X_i,X_j,X_k}
    F(X_i,X_j,X_k)=g\left([X_i,JX_j]+[X_j,JX_i],X_k\right)\\[4pt]
    \phantom{\mathop{\makebox{\huge $\sigma$}} \limits_{X_i,X_j,X_k}
    F(X_j,X_j,X_k)=}+g\left([X_j,JX_k]+[X_k,JX_j],X_i\right)\\[4pt]
    \phantom{\mathop{\makebox{\huge $\sigma$}} \limits_{X_i,X_j,X_k}
    F(X_j,X_j,X_k)=}+g\left([X_k,JX_i]+[X_i,JX_k],X_j\right).
\end{array}
\end{equation}

The condition for $(G,J,g)$ to be of $\W_3$-manifold is the
vanishing of \eqref{F_3}. Having in mind \eqref{J}, \eqref{g} and
the skew-symmetry of the Lie brackets, we determine at the first
stage the constants $C_{ij}^k$ from the vanishing of the right
side of \eqref{F_3} and hence we obtain the commutators of the
basis vector fields as follows
\begin{equation}\label{[]}
\begin{array}{l}
[X_1,X_3]=\lambda_2 X_2+\lambda_4 X_4,\\[4pt]
[X_2,X_4]=\lambda_1 X_1+\lambda_3 X_3,\\[4pt]
[X_2,X_3]=\lambda_5 X_1+\lambda_6 X_2+\lambda_7 X_3+\lambda_8 X_4,\\[4pt]
[X_3,X_4]=\lambda_9 X_1+\lambda_{10} X_2+\lambda_{11} X_3+\lambda_{12} X_4,\\[4pt]
[X_4,X_1]=(\lambda_2+\lambda_5) X_1+(\lambda_1+\lambda_6) X_2
+(\lambda_4+\lambda_7) X_3+(\lambda_3+\lambda_8) X_4,\\[4pt]
[X_2,X_1]=(\lambda_9+\lambda_4) X_1+(\lambda_{10}-\lambda_3) X_2
+(\lambda_{11}-\lambda_2) X_3+(\lambda_{12}+\lambda_1) X_4,\\[4pt]
\end{array}
\end{equation}
where the coefficients $\lambda_i\in \R$ $(i=1,2,\dots,12)$ must
satisfy the Jacobi condition \eqref{Jac}.

Conversely, let us consider an almost complex manifold with Norden
metric $(G,J,g)$ defined by the above manner, where the conditions
\eqref{[]} hold. We check immediately that the right side of
\eqref{F_3} vanishes. Then the left side of \eqref{F_3} is zero
and according \eqref{class} we obtain that $(G,J,g)$ belongs to
the class $\W_3$.

Now we put the condition the Norden metric $g$ to be a Killing
metric of the Lie group $G$ with the corresponding Lie algebra
$\g$, i.e.
\begin{equation}\label{Kil}
    g(\ad X(Y),Z)=-g(Y,\ad X(Z),
\end{equation}
where $X,Y,Z \in \g$ and $\ad X (Y)=[X,Y]$. It is equivalent to
the condition the metric $g$ to be an invariant metric, i.e.
\begin{equation}\label{inv}
    g\left([X,Y],Z\right)+g\left([X,Z],Y\right)=0.
\end{equation}

Then, according to \eqref{inv}, the equations \eqref{[]} take the
following form
\begin{equation}\label{[]4}
\begin{array}{l}
[X_1,X_3]=\lambda_2 X_2+\lambda_4 X_4,\\[4pt]
[X_2,X_4]=\lambda_1 X_1+\lambda_3 X_3,\\[4pt]
[X_2,X_3]=-\lambda_2 X_1-\lambda_3 X_4,\\[4pt]
[X_3,X_4]=-\lambda_4 X_1+\lambda_{3} X_2,\\[4pt]
[X_4,X_1]=\lambda_1 X_2+\lambda_4 X_3,\\[4pt]
[X_2,X_1]=-\lambda_2 X_3+\lambda_1 X_4,\\[4pt]
\end{array}
\end{equation}

By direct verification we prove that the commutators from
\eqref{[]4} satisfy the Jacobi identity. The Lie groups $G$ thus
obtained are of a family which is characterized by four parameters
$\lambda_i$ $(i = 1,\dots, 4)$.

Therefore, for the manifold $(G,J,g)$ constructed above, we
establish the truthfulness of the following
\begin{thm}
Let $(G,J,g)$ be a 4-dimensional almost complex manifold with
Norden metric, where $G$ is a connected Lie group with
corresponding Lie algebra $\g$ determined by the global basis of
left invariant vector fields $\{X_1,X_2,X_3,X_4\}$; $J$ is an
almost complex structure defined by \eqref{J} and $g$ is an
invariant Norden metric determined by \eqref{g} and \eqref{inv}.
Then $(G,J,g)$ is a quasi-K\"ahler manifold with Norden metric if
and only if $G$ belongs to the 4-parametric family of Lie groups
determined by the conditions \eqref{[]4}.
\end{thm}


\section{Geometric characteristics of the constructed manifold}

Let $(G,J,g)$ be the 4-dimensional quasi-K\"ahler manifold with
Norden metric introduced in the previous section.

From the Levi-Civita connection \eqref{LC} and the condition
\eqref{inv} we have
\begin{equation}\label{invLC}
    \nabla_{X_i} X_j=\frac{1}{2}[X_i,X_j]\quad (i,j=1,2,3,4).
\end{equation}

\subsection{The components of the tensor F}
Then by direct calculations, having in mind \eqref{F}, \eqref{J},
\eqref{g}, \eqref{[]4} and \eqref{invLC}, we obtain the nonzero
components of the tensor $F$ as follows

\begin{equation}\label{Fijk}
\begin{split}
-F_{122}&=-F_{144}=2F_{212}=2F_{221}=2F_{234}\\[4pt]
\phantom{-F_{122}}
&=2F_{243}=2F_{414}=-2F_{423}=-2F_{432}=2F_{441}=\lambda_1,\\[4pt]
2F_{112}&=2F_{121}=2F_{134}=2F_{143}=-2F_{211}\\[4pt]
\phantom{2F_{112}}
&=-2F_{233}=-2F_{314}=2F_{323}=2F_{332}=-2F_{341}=\lambda_2,\\[4pt]
2F_{214}&=-2F_{223}=-2F_{232}=2F_{241}=F_{322}\\[4pt]
\phantom{2F_{112}}
&=F_{344}=-2F_{412}=-2F_{421}=-2F_{434}=-2F_{443}=\lambda_3,\\[4pt]
-2F_{114}&=2F_{123}=2F_{132}=-2F_{141}=-2F_{312}\\[4pt]
\phantom{-2F_{114}}
&=-2F_{321}=-2F_{334}=-2F_{343}=F_{411}=F_{433}=\lambda_4,\\[4pt]
\end{split}
\end{equation}
where $F_{ijk}=F(X_i,X_j,X_k)$.

\subsection{The square norm of the Nijenhuis
tensor}

Let $N$ be the Nijenhuis tensor of the almost complex structure
$J$ on $G$, i.e.
\begin{equation}\label{N}
    N(X,Y)=[X,Y]+J[JX,Y]+J[X,JY]-[JX,JY], \quad X, Y \in \g.
\end{equation}
Having in mind \eqref{[]4} we obtain the nonzero components
$N_{ij}=N(X_i,X_j)$ $(i, j=1,2,3,4)$ as follows
\begin{equation}\label{Nij}
\begin{array}{l}
    N_{12}=-N_{34}=2\left(\lambda_4 X_1-\lambda_3 X_2+\lambda_2 X_3-\lambda_1
        X_4\right),\\[4pt]
    N_{14}=-N_{23}=2\left(\lambda_2 X_1-\lambda_1 X_2-\lambda_4
        X_3+\lambda_3 X_4\right).\\[4pt]
\end{array}
\end{equation}
Therefore its square norm $\nN=g^{ik}g^{ks}g(N_{ij},N_{ks}$ has
the form
\begin{equation}\label{nN}
    \nN=-32\left(\lambda_1^2+\lambda_2^2-\lambda_3^2-\lambda_4^2\right),
\end{equation}
where the inverse matrix of $g$ has the form
\begin{equation}\label{g^ij}
    \left(g^{ij}\right)=\left(%
\begin{array}{cccc}
  1 & 0 & 0 & 0 \\
  0 & 1 & 0 & 0 \\
  0 & 0 & -1 & 0 \\
  0 & 0 & 0 & -1 \\
\end{array}%
\right)
\end{equation}

\begin{prop}\label{Prop3.1}
    The manifold $(G,J,g)$ is isotropic K\"ahlerian if and only if
    its Nijenhuis tensor is isotropic.
\end{prop}

\subsection{The square norm of $\nabla J$} According to \eqref{g} and \eqref{g^ij},
from \eqref{snorm} we obtain the square norm of $\nabla J$ as
\begin{equation}\label{snorm4}
    \nJ=4\left(\lambda_1^2+\lambda_2^2-\lambda_3^2-\lambda_4^2\right).
\end{equation}

The last equation in accordance with \dfnref{iK} implies the
following
\begin{prop}\label{Th1}
    The manifold $(G,J,g)$ is isotropic K\"ahlerian if and only if
    the condition $\lambda_1^2+\lambda_2^2-\lambda_3^2-\lambda_4^2=0$
    holds.
\end{prop}

\begin{rem}
The last theorem means that the set of vectors with the
coordinates $(\lambda_1,\lambda_2,\lambda_3,\lambda_4)$ at an
arbitrary point $p\in G$ describes the isotropic cone in $T_pG$.
\end{rem}

\subsection{The components of $R$}
Let $R$ be the curvature tensor of type (0,4) determined by
\eqref{R04} and \eqref{R} on $(G,J,g)$. We denote its components
by $R_{ijks}=R(X_i,X_j,X_k,X_s)$ $(i,j,k,s=1,2,3,4)$. Using
\eqref{invLC} and the Jacobi identity we receive
\begin{equation}\label{invRijks}
    R_{ijks}=-\frac{1}{4}g\Bigl(\bigl[[X_i,X_j],X_k\bigr],X_s\Bigr).
\end{equation}
From the last formula and \eqref{[]4} we get the nonzero
components of $R$ as follows
\begin{equation}\label{Rijks}
\begin{array}{ll}
    R_{1221}=-\frac{1}{4}\left(\lambda_1^2+\lambda_2^2\right),\quad
    &
    R_{1331}=\frac{1}{4}\left(\lambda_2^2-\lambda_4^2\right),\\[4pt]
    R_{1441}=-\frac{1}{4}\left(\lambda_1^2-\lambda_4^2\right),\quad
    &
    R_{2332}=\frac{1}{4}\left(\lambda_2^2-\lambda_3^2\right),\\[4pt]
    R_{2442}=\frac{1}{4}\left(\lambda_1^2-\lambda_3^2\right),\quad
    &
    R_{3443}=\frac{1}{4}\left(\lambda_3^2+\lambda_4^2\right),\\[4pt]
    R_{1341}=R_{2342}=-\frac{1}{4}\lambda_1\lambda_2,\quad
    &
    R_{2132}=-R_{4134}=\frac{1}{4}\lambda_1\lambda_3,\\[4pt]
    R_{1231}=-R_{4234}=\frac{1}{4}\lambda_1\lambda_4,\quad
    &
    R_{2142}=-R_{3143}=\frac{1}{4}\lambda_2\lambda_3,\\[4pt]
    R_{1241}=-R_{3243}=\frac{1}{4}\lambda_2\lambda_4,\quad
    &
    R_{3123}=R_{4124}=\frac{1}{4}\lambda_3\lambda_4.\\[4pt]
\end{array}
\end{equation}

\subsection{The components of $\rho$ and the value of $\tau$}
Having in mind \eqref{rho-tau} and \eqref{g^ij} we obtain the
components $\rho_{ij}=\rho(X_i,X_j)$ $(i,j=1,2,3,4)$ of the Ricci
tensor and the scalar curvature $\tau$
\begin{equation}\label{rho_ij}
\begin{array}{c}
\begin{array}{ll}
    \rho_{11}=-\frac{1}{2}\left(\lambda_1^2+\lambda_2^2-\lambda_4^2\right),\quad
    &
    \rho_{22}=-\frac{1}{2}\left(\lambda_1^2+\lambda_2^2-\lambda_3^2\right),\\[4pt]
    \rho_{33}=\frac{1}{2}\left(\lambda_2^2-\lambda_3^2-\lambda_4^2\right),\quad
    &
    \rho_{44}=\frac{1}{2}\left(\lambda_1^2+\lambda_3^2-\lambda_4^2\right),\\[4pt]
\end{array}
\\[4pt]
\begin{array}{lll}
    \rho_{12}=-\frac{1}{2}\lambda_3\lambda_4,\quad
    &
    \rho_{13}=\frac{1}{2}\lambda_1\lambda_3,\quad
    &
    \rho_{14}=\frac{1}{2}\lambda_2\lambda_3,\\[4pt]
    \rho_{23}=\frac{1}{2}\lambda_1\lambda_4,\quad
    &
    \rho_{24}=\frac{1}{2}\lambda_2\lambda_4,\quad
    &
    \rho_{34}=-\frac{1}{2}\lambda_1\lambda_2;\\[4pt]
\end{array}
\end{array}
\end{equation}
\begin{equation}\label{tau}
    \tau=-\frac{3}{2}\left(\lambda_1^2+\lambda_2^2-\lambda_3^2-\lambda_4^2\right),
\end{equation}
i.e. the scalar curvature on $(G,J,g)$ is constant.

The last equation and \propref{Th1} imply immediately
\begin{prop}\label{Th2}
    The manifold $(G,J,g)$ is isotropic K\"ahlerian if and only if it is scalar flat.
\end{prop}

\subsection{The Weyl tensor}

Now, let us consider the Weyl tensor $W$ on $(G,J,g)$ defined by
\eqref{W} and \eqref{psi-pi}. Taking into account \eqref{g},
\eqref{Rijks}, \eqref{rho_ij} and \eqref{tau}, we establish the
truthfulness of the following
\begin{thm}\label{Th3}
The manifold $(G,J,g)$ has vanishing Weyl tensor.
\end{thm}
The \thmref{Prop1.1} and \thmref{Th3} imply immediately the last
\begin{thm}\label{Th4}
The manifold $(G,J,g)$ is locally symmetric.
\end{thm}

\subsection{The sectional curvatures and the holomorphic bisectional curvature}

Let us consider the characteristic 2-planes $\alpha_{ij}$ spanned
by the basis vectors $\{X_i,X_j\}$ at an arbitrary point of the
manifold:
\begin{itemize}
    \item holomorphic 2-planes - $\alpha_{13}$, $\alpha_{24}$;
    \item totally real 2-planes - $\alpha_{12}$, $\alpha_{14}$, $\alpha_{23}$, $\alpha_{34}$.
\end{itemize}
Then, using \eqref{psi-pi}, \eqref{k}, \eqref{g} and
\eqref{Rijks}, we obtain the corresponding sectional curvatures
\begin{equation}\label{k_ik}
    \begin{array}{ll}
    k(\alpha_{13})=-\frac{1}{4}\left(\lambda_2^2-\lambda_4^2\right),\quad
    &
    k(\alpha_{24})=-\frac{1}{4}\left(\lambda_1^2-\lambda_3^2\right),\\[4pt]
    k(\alpha_{12})=-\frac{1}{4}\left(\lambda_1^2+\lambda_2^2\right),\quad
    &
    k(\alpha_{14})=-\frac{1}{4}\left(\lambda_1^2-\lambda_4^2\right),\\[4pt]
    k(\alpha_{23})=-\frac{1}{4}\left(\lambda_2^2-\lambda_3^2\right),\quad
    &
    k(\alpha_{14})=\frac{1}{4}\left(\lambda_3^2+\lambda_4^2\right).\\[4pt]
\end{array}
\end{equation}

Taking into account \eqref{h}, \eqref{g} and \eqref{Rijks}, we
obtain that the holomorphic bisectional curvature of the unique
pair of basis holomorphic 2-planes $\{\alpha_{13}, \alpha_{24}\}$
vanishes, i.e.
\begin{equation}\label{h=0}
    h(X_1,X_2)=0.
\end{equation}

\subsection{The isotropic-K\"ahlerian property} Having in mind the Propositions \ref{Prop3.1}--\ref{Th2}  in the previous
subsections, we give the following
\begin{thm}
    The following conditions are equivalent for the manifold $(G,J,g)$:
    \begin{enumerate}
    \renewcommand{\labelenumi}{(\roman{enumi})}
    \item
    $(G,J,g)$ is isotropic K\"ahler manifold;
    \item
    the Nijenhuis tensor is isotropic;
    \item
    the condition $\lambda_1^2+\lambda_2^2-\lambda_3^2-\lambda_4^2=0$
    holds;
    \item
    the scalar curvature vanishes.
    \end{enumerate}
\end{thm}


\end{document}